\documentclass[11pt,a4paper]{amsart}
\usepackage{epsfig,a4wide}

\newtheorem{theorem}{Theorem}
\newtheorem{lemma}[theorem]{Lemma}

\begin{document}

\title{Dehn fillings creating essential spheres and tori}

\author[S. Lee]{Sangyop Lee}
\address{Department of Mathematics, Korea Advanced Institute of
         Science and Technology, Taejon 305--701, Korea}
\email{apple@mathx.kaist.ac.kr}
\author[S. Oh]{Seungsang Oh}
\address{Department of Mathematics, Chonbuk National University,
         Chonju, Chonbuk 561-756, Korea}
\email{soh@math.chonbuk.ac.kr}
\author[M. Teragaito]{Masakazu Teragaito}
\address{Department of Mathematics and Mathematics Education,
         Hiroshima University, Kagamiyama 1-1-1,
         Higashi-Hiroshima 739-8524, Japan}
\email{teragai@hiroshima-u.ac.jp}

\begin{abstract}
Let $M$ be a simple $3$-manifold with a toral boundary component. 
It is known that if two Dehn fillings on $M$ along 
the boundary produce a reducible manifold and a toroidal manifold, 
then the distance between the filling slopes is at most three.
This paper gives a remarkably short proof of this result. 
\end{abstract}

\maketitle

Let $M$ be a {\em simple\/} $3$-manifold with a toral boundary 
component $\partial_0 M$, 
that is, it contains no essential sphere, disk, torus or annulus. 
The {\em slope\/} of an essential, unoriented, simple closed curve 
on $\partial_0 M$ is its isotopy class. 
We assume that $\alpha$ and $\beta$ are two slopes on $\partial_0 M$ 
such that $M(\alpha)$ is a reducible manifold and 
$M(\beta)$ contains an essential torus.
The goal is to measure the upper bound for $\Delta(\alpha,\beta)$ 
(the minimal geometric intersection number of $\alpha$ and $\beta$),
and  Oh \cite{Oh1} and Wu \cite{Wu} independently gave the 
optimum upper bound 3 for this case.
In the present paper we give another short proof of this result.

\begin{theorem}\label{thm:main}
Let $M$ be a simple $3$-manifold such that $M(\alpha)$ is
reducible and $M(\beta)$ is toroidal. Then $\Delta(\alpha,\beta)\le 3$.
\end{theorem}

Assume $\Delta(\alpha,\beta) \ge 4$ for contradiction. Let $\widehat{Q}$
be a reducing sphere in $M(\alpha)$ which intersects the filling
solid torus $V_{\alpha}$ in a family of meridian disks. 
We choose $\widehat{Q}$ so that $q=|\widehat{Q}\cap V_{\alpha}|$ is minimal 
(over all reducing spheres $\widehat{Q}$ in $M(\alpha)$). 
Similarly let $\widehat{T}$ be an essential torus in $M(\beta)$ 
meeting in $t=|\widehat{T}\cap V_{\beta}|$ meridian disks, 
the number of which is minimal over all such tori.
Let $Q=\widehat{Q}\cap M$ and $T=\widehat{T}\cap M$. By an isotopy of $Q$, 
we may assume that $Q$ and $T$ intersect transversely, and
$Q \cap T$ has the minimal number of components. 
Then as described in \cite{Oh1}, 
we obtain graphs $G_Q$ in $\widehat{Q}$ and $G_T$ in $\widehat{T}$. 
We use the definitions and terminology of \cite[Section2]{Oh1}.

\begin{lemma} \cite[Lemmas 2.3, 4.1, 4.3]{BZ} \
$q \ge 3$ and $t \ge 3$.
\end{lemma}

\begin{lemma}\label{lem:Schar} \cite[Theorem 2.4]{GL2} \
$G_T$ cannot contain  Scharlemann cycles on distinct label pairs.
\end{lemma}

\begin{lemma}\label{lem:parallel} \cite[Lemmas 1.4, 1.8]{Wu} \
$G_Q$ cannot have more than $\frac{t}{2}+2$ mutually parallel edges
connecting parallel vertices. 
Furthermore, if $G_Q$ has $\frac{t}{2}+2$ mutually parallel edges 
connecting parallel vertices, then $t \equiv 0 \ (\mbox{mod} \ 4)$.
\end{lemma} 

\begin{theorem}\label{lem:JLOT} \cite[Theorem 1.1]{JLOT}
If $M(\alpha)$ and $M(\beta)$ contain a projective plane  and an
essential torus respectively, then either $\Delta(\alpha,\beta)\le 2$
or $\Delta(\alpha,\beta)=3$ with $t=2$.
\end{theorem}  

Now we will prove Theorem \ref{thm:main}.

\begin{proof}
We distinguish two cases;
\vspace{2mm}

\noindent {\bf Case I)} Suppose that there is a vertex $x$ of $G_Q$ such that
more than $\frac{3}{2}t$  edges connect $x$ to antiparallel vertices.

This implies that there exist more than $\frac32t$ $x$-edges 
connecting parallel vertices in $G_T$. 
Consider the subgraph of $G_T$ consisting all vertices and 
those all $x$-edges of $G_T$.
An Euler characteristic count gives that it contains more than $\frac12t$ disk
faces whose boundaries are $x$-edge cycles of $G_T$. 
By \cite[Proposition 5.1]{HM}, each disk face contains a Scharlemann cycle. 
By Lemma \ref{lem:Schar},
these all Scharlemann cycles are, say, $12$-Scharlemann cycles.

Construct a graph $\Gamma$ in $\widehat{T}$ as follows. 
Choose a {\em dual vertex\/} in the interior of the face $G_T$ bounded 
by a $12$-Scharlemann cycle, 
and let the vertices of $\Gamma$ be the vertices of $G_T$ together 
with these  dual vertices.
The edges of $\Gamma$ are defined by joining each dual vertex to 
the vertices of the corresponding Scharlemann cycle in the obvious way. 
Let $s$ be the number of $12$-Scharlemann cycles. 
Then $s> \frac12t$. 
$\Gamma$ has $t+s$ vertices and at least $3s$ edges 
because each Scharlemann cycle has order  at least 3 by Theorem \ref{lem:JLOT}. 
Note that if $G_T$ contains an Scharlemann cycle of length two, 
then $M(\alpha)$ contains a projective plane \cite{GLi}.
Again an Euler characteristic count guarantees that
$\Gamma$ has a disk face $E$. But $E$ determines a $1$-edge cycle bounding
a disk face in $E$ which, as long as $q>2$, contains a Scharlemann cycle. 
This contradicts the definition of $\Gamma$.
\vspace{2mm}
 
\noindent {\bf Case II)} As the negation of case I), 
suppose that each vertex $x$ of $G_Q$ has at least $(\Delta-\frac32)t$ labels 
where edges connecting parallel vertices are incident.

Let $G_Q^+$be the subgraph of $G_Q$ consisting of all vertices and edges
connecting parallel vertices of $G_Q$.
By the assumption, every vertex of $G_Q^+$ has valency
at least $(\Delta-\frac32)t$.
Let $\overline G_Q^+$ be the reduced graph of $G_Q^+$.
By Lemma \ref{lem:parallel}, any vertex has valency at least $3$ in 
$\overline G_Q^+$.
So, we can choose an innermost component $H$ of $\overline G_Q^+$, and
a block $\overline \Lambda$ of $H$ with at most one cut vertex of $H$.
Let $\Lambda$ be the subgraph of $G_Q^+$ corresponding to $\overline \Lambda$.  

Note that there is a disk $D$ in $\widehat{Q}$ such that $D\cap G_Q^+=\Lambda$.
A vertex of $\overline{\Lambda}$ is a \textit{boundary vertex\/} if
there is an arc connecting it to $\partial D$ whose interior is disjoint from
$\overline{\Lambda}$,
and an \textit{interior vertex\/} otherwise.
Remark that $\overline{\Lambda}$ has an interior vertex by \cite[Lemma 3.2]{Wu}.
 
If $\Lambda$ contains a cut vertex of $H$, then let $x$ 
be a label such that the cut vertex
has an edge attached with label $x$ there.
Otherwise, let $x$ be any label.
Then each vertex of $\Lambda$, except a cut vertex, has at least
two edges attached with label $x$.
Thus we can find a great $x$-cycle.
This implies that $\Lambda$ has a Scharlemann cycle
and so $t$ is even by \cite[Lemma 2.2(1)]{BZ}.

Let $v,e$ and $f$ be the numbers of vertices, 
edges and faces of its reduced graph $\overline{\Lambda}$. 
Also say $v_i,v_{\partial}$ and $v_c$ the numbers of interior vertices, 
boundary vertices and a cut vertex of $\overline{\Lambda}$. 
Hence $v=v_i+v_{\partial}$ and $v_c=0$ or $1$. 
Since each face of $\overline{\Lambda}$ is a disk with at least $3$ sides, 
we have $2e\ge3f+v_{\partial}$. 
Combined with $1=\chi(\mathrm{disk})= v-e+f$, 
we get $e\le3v-v_{\partial}-3=3v_i+2v_{\partial}-3$.

Suppose that every interior vertex of $\overline{\Lambda}$
has valency at least $6$ and that every boundary vertex except a cut vertex
has valency at least $4$. 
Since a cut vertex has valency at least $2$, 
we have $2e\ge 6v_i+4(v_{\partial}-v_c)+2v_c$. 
These two inequalities give us that $0\le v_c-3$, a contradiction.

Hence there are two cases. 
For the first case, assume that some  boundary vertex $y$ of 
$\overline{\Lambda}$ has valency at most $3$.
There are at least $(\Delta -\frac32)t$ edges in $G_Q$ 
which are incident to $y$ and connect $y$ to parallel vertices. 
By Lemma \ref{lem:parallel}, 
$\frac52t \le (\Delta -\frac32) t \le 3(\frac{t}{2}+2)$, 
i.e. $t \leq 6$ and $t \neq 6$. 

But in the case of $t=4$ there are at least $4$ mutually 
parallel edges in $\Lambda$,
which implies that $G_Q$ contains two S-cycles with disjoint 
label pairs.
Then $M(\beta)$ contains a Klein bottle
as in the proof of \cite[Lemma 3.10]{GL1}.
(Note that $M(\beta)$ is irreducible \cite{GL2}.
Hence the edges of an $S$-cycle cannot lie in a disk on $\widehat{T}$.)
This contradicts to \cite[Theorem 1.1]{Oh2}. 

For the second case, assume that some interior vertex of 
$\overline{\Lambda}$ has valency at most $5$.
Again we have $4t \le \Delta t \le 5(\frac{t}{2}+2)$.
Then we can use the same argument as above.
\end{proof}


\begin{thebibliography}{CGLS}
\bibitem[BZ]{BZ} S. Boyer and X. Zhang,
            {\em Reducing Dehn filling and toroidal Dehn filling},
            Topology Appl. \textbf{68} (1996), 285--303.
\bibitem[CGLS]{CGLS} M. Culler, C. Gordon, J. Luecke and P. Shalen,
            {\em Dehn surgery on knots},
            Ann. Math. \textbf{125} (1987), 237--300.
\bibitem[GL1]{GL1} C. Gordon and J. Luecke,
            {\em Dehn surgeries on knots creating essential tori, I},
            Communications in Analysis and Geometry
            \textbf{3} (1995), 597--644.
\bibitem[GL2]{GL2} C. Gordon and J. Luecke,
            {\em Reducible manifolds and Dehn surgery},
            Topology \textbf{35} (1996), 385--409.
\bibitem[GLi]{GLi} C. Gordon and R. Litherland,
            {\em Incompressible planar surfaces in 3-manifolds},
            Topology Appl. \textbf{18} (1984), 121--144.
\bibitem[HM]{HM} C. Hayashi and K. Motegi,
            {\em Only single twists on unknots can produce
            composite knots},
            Trans. Amer. Math. Soc. \textbf{349} (1997), 4465--4479.
\bibitem[JLOT]{JLOT} G. T. Jin, S. Lee, S. Oh and M. Teragaito,
            {\em $P^2$-reducing and toroidal Dehn fillings},
            preprint.
\bibitem[Oh1]{Oh1} S. Oh, {\em Reducible and toroidal
            3-manifolds obtained by Dehn fillings},
            Topology Appl. \textbf{75} (1997), 93--104.
\bibitem[Oh2]{Oh2} S. Oh, {\em Dehn filling, reducible
            3-manifolds, and Klein bottles},
            Proc. Amer. Math. Soc. \textbf{126-1} (1998), 289--296.
\bibitem[Wu]{Wu} Y.Q. Wu, {\em Dehn fillings producing
            reducible manifolds and toroidal manifolds},
            Topology \textbf{37} (1998), 95--108.
\end{thebibliography}
\end{document}